\documentclass[12pt,reqno]{amsart}

\usepackage[final%,
]{showkeys}
\usepackage{AKstyle}

\usepackage[table]{xcolor}

\newcommand{\vertiii}[1]{{\left\vert\kern-0.25ex\left\vert\kern-0.25ex\left\vert #1
    \right\vert\kern-0.25ex\right\vert\kern-0.25ex\right\vert}}

\numberwithin{equation}{section}

\usepackage{caption}
\usepackage[labelfont=rm]{subcaption}

\usepackage[pagebackref=true, colorlinks]{hyperref}

\hypersetup{pdffitwindow=true,linkcolor=blue,citecolor=blue,urlcolor=blue,menucolor=blue}

\usepackage{comment}

\long\def\ignore#1{}

\begin{document}%%

\author{Arthur Baragar}
%\thanks{Baragar is partially supported by ...}
\email{baragar@unlv.nevada.edu}
\address{Department of Mathematical Sciences, University of Nevada Las Vegas, Las Vegas, NV 89154}
\author{Alex Kontorovich}
\thanks{Kontorovich is partially supported by
an NSF CAREER grant DMS-1455705, an NSF FRG grant DMS-1463940, and a BSF grant.}
\email{alex.kontorovich@rutgers.edu}
\address{Department of Mathematics, Rutgers University, New Brunswick, NJ 08854}

\title
%{An efficient construction of four mutually tangent circles}
{Efficiently Constructing %Mutually
Tangent Circles}

\date{\today}
\maketitle

\section{Introduction}

The famous Problem of Apollonius is to construct a circle tangent to three given ones in a plane. The three circles may also be limits of circles, that is, points or lines, and ``construct'' of course refers to straightedge and compass. In this note, we consider the problem of constructing tangent circles from the point of view of {\it efficiency}. By this we mean using as few {\it moves} as possible, where a move is the act of drawing a line or circle. (Points are free as they do not harm the straightedge or compass, and all lines are considered endless, so there is no cost to ``extending'' a line segment.)
Our goal is to present, in what we believe is the most efficient way possible, a construction of four mutually tangent circles.
(Five circles of course cannot be mutually tangent in the plane, for their tangency graph, the complete graph $K_{5}$, is non-planar.)
We first present our construction before giving some remarks comparing it to others we found in the literature.

\section{Baby Cases: One and Two Circles}

\begin{figure}
\includegraphics[scale=.94]{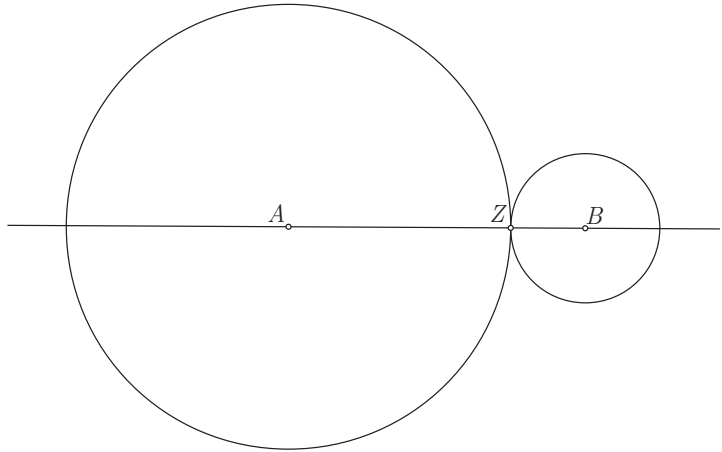}
\caption{Two tangent circles.}
\label{fig:1}
\end{figure}

Constructing one circle obviously costs one move: let $A$ and $Z$ be any distinct points in the plane and draw the circle $O_{A}$ with center $A$ % and radius $|AZ|$.
and passing through $Z$.
Given $O_{A}$, constructing a second circle tangent to it costs two more moves: draw a line through $AZ$, % (tally: 1 move),
and put an arbitrary point $B$ on this line (say, outside $O_{A}$). Now draw the circle $O_{B}$
%(tally: 2 moves)
with center $B$ and %radius $|BZ|$
passing through $Z$; then $O_{A}$ and $O_{B}$ are obviously tangent at $Z$, see \figref{fig:1}.
It should be clear that one cannot do better than two moves, for otherwise one could draw the circle $O_{B}$ immediately; but this  requires knowledge of a point on $O_{B}$.

\section{Warmup: Three Circles}

Given \figref{fig:1}, that is, the two circles $O_{A}$ and $O_{B}$, tangent at $Z$, and the line $AB$, how many moves does it take to construct a third  circle tangent to both $O_{A}$ and $O_{B}$? We encourage %the
readers at this point to stop and try this problem
themselves.
% herself.

\pagebreak

\begin{prop}
Given \figref{fig:1}, a circle tangent to both $O_{A}$ and $O_{B}$ is constructible in at most \emph{five} moves.
\end{prop}
%Our solution requires five moves:

We first give the construction, then the proof that it works.

\subsection*{The Construction}
\begin{figure}
\includegraphics[scale=.94]{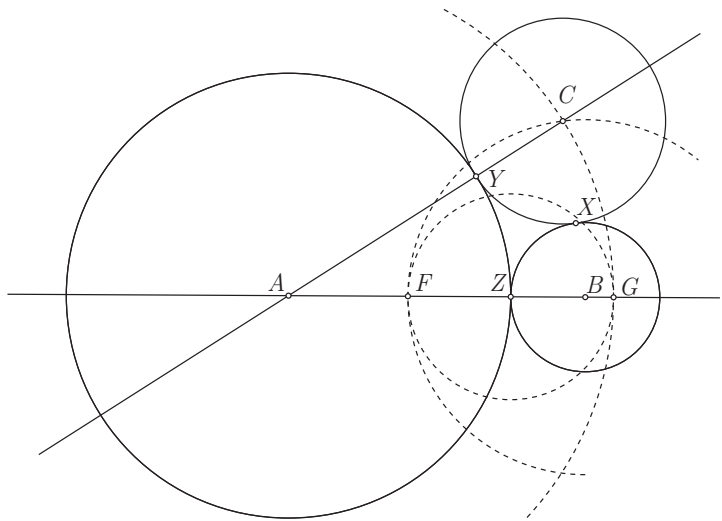}
\caption{A third tangent circle.  The solid lines/circles are the initial and final objects, while the dotted figures are the intermediate constructions. The point $X$ lies on the line $BC$ (not shown) and is the point of tangency between $O_B$ and $O_C$.  }
\label{fig:2}
\end{figure}

Draw an arbitrary circle $O_{Z}$ centered at $Z$  (this is move 1), and let it intersect $AB$ at $F$ and $G$, say, with $A$ and $F$ on the same side of $Z$.  Next draw the circle  centered at $A$ and passing through $G$ (move 2), and the circle centered at $B$ through $F$ (move 3); see \figref{fig:2}.  Let these two circles intersect at $C$.  Construct the line $AC$ (move 4) and let it intersect $O_A$ at $Y$.  Finally, draw the circle $O_{C}$ centered at $C$ %with radius $|CY|$
and passing through $Y$
(move 5); then $O_{C}$ is tangent to  $O_B$ at $X$, say.

\subsection*{The Proof}
%Let $r=|ZF|$ be the radius of $O_Z$.  It is clear that $O_C$ is tangent to $O_A$ at $Y$, and by construction, has radius $r$.  Let $BC$ (not constructed) intersect $O_B$ at $X$.  Then $|CX|=|BC|-|BX|=|BF|-|BZ|=r$, so $X$ lies on $O_C$.  Hence, $O_C$ is tangent to $O_B$ at $X$.
It is elementary to verify that the above construction works, and that the radius of $O_{C}$ is the same as that of $O_{Z}$.
Note in fact that the locus of all centers $C$ of circles $O_C$ tangent to both $O_{A}$ and $O_{B}$ forms a hyperbola with foci $A$ and $B$. Indeed, let the circles $O_{A}$, $O_{B}$, and $O_{C}$ have radii $a$, $b$, and $c$, resp.; then $|AC| =a+c$ and $|BC| =b+c$, so $|AC|-|BC| =a-b$ is constant for any choice of $c$.
\hfill $\sq$

\ignore{We claim that the locus of all centers $C$ of circles $O_{C}$ tangent to both $O_{A}$ and $O_{B}$ forms a hyperbola with foci $A$ and $B$.
Indeed, let the circles $O_{A}$, $O_{B}$, and $O_{C}$ have radii $a$, $b$, and $c$, resp.; then $AC=a+c$ and $BC=b+c$, so $AC-BC=a-b$ is fixed. It is then elementary to verify that the above construction works. % to produce $C$.
}

\section{Main Theorem: the Fourth Circle}

Finally we come to the main event, the fourth tangent circle, which we call the Apollonian circle.\footnote{Many objects in the literature are named after Apollonius though he had nothing to do with them, such as the Apollonian gasket and the Apollonian group (see, e.g., \cite{Kontorovich2013}). The fourth tangent circle really {\it is} due to him, though most authors refer to it as the ``Soddy'' circle.}
We are given three mutually tangent circles, $O_{A}$, $O_{B}$ and $O_{C}$, lines $AB$ and $AC$, and the points of tangency $X$, $Y$, and $Z$; that is, we are given the already constructed objects in \figref{fig:2}.
% (though we will not need the intermediate constructions other than $AC$).

\begin{thm}
An Apollonian circle tangent to $O_{A}$, $O_{B}$ and $O_{C}$ in \figref{fig:2} is constructible in at most \emph{seven} %$elementary
moves.
\end{thm}

\subsection*{The construction}  
\begin{figure}
\includegraphics[scale=.94]{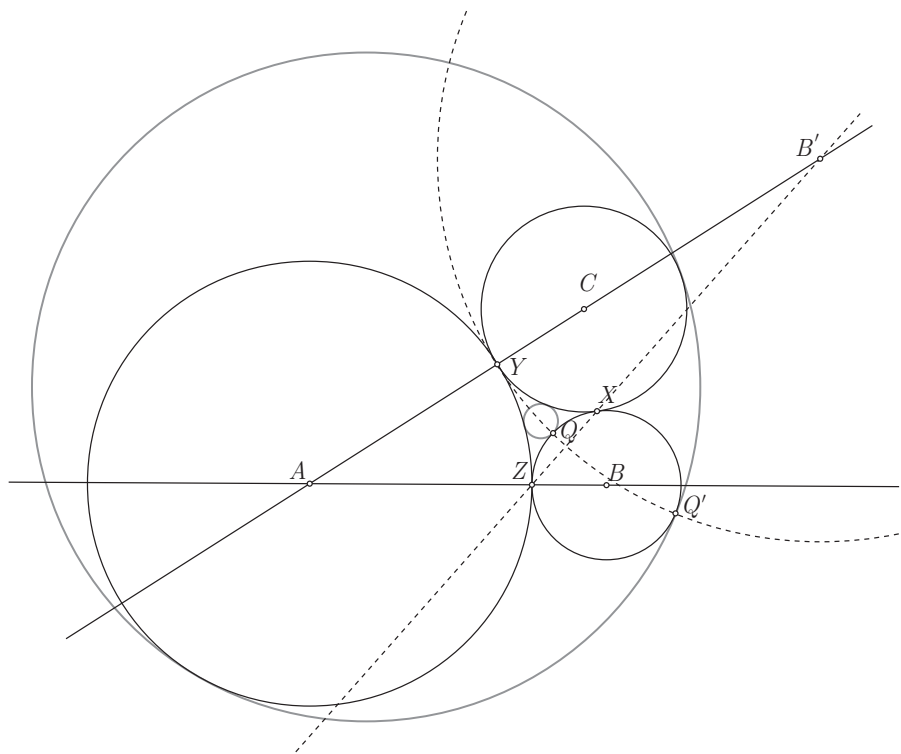}
\caption{The circle $O_{B'}$ and points $Q$ and $Q'$.  The grey circles are the Apollonian circles $O_S$ and $O_{S'}$ that we are in the process of constructing.}
\label{fig:3}
\end{figure}

\begin{figure}
\includegraphics[scale=.94]{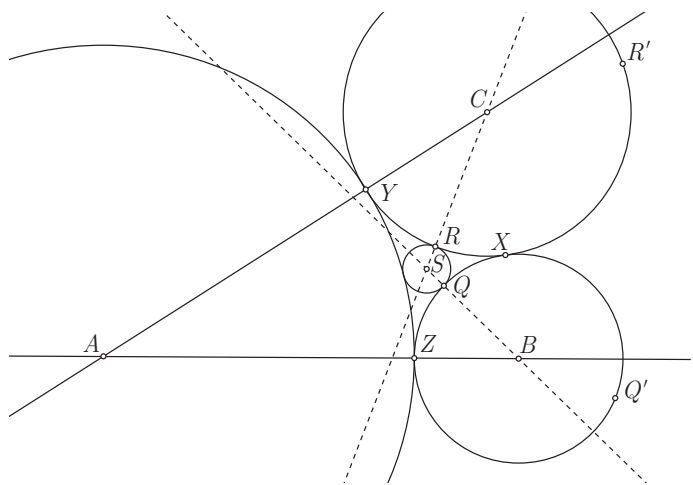}
\caption{The construction of $S$ and the Apollonian circle $O_{S}$.}
\label{fig:4}
\end{figure}

Draw the line $XZ$ (this is move 1) and let it intersect $AC$ at $B'$.  Draw the circle $O_{B'}$   centered at $B'$ %with radius $|B'Y|$ (%this is
and passing through $Y$
(move 2).  It intersects $O_B$ at $Q$ and $Q'$; see \figref{fig:3}. We repeat this procedure: draw the line $XY$, let it intersect $AB$ at $C'$, draw the circle $O_{C'}$ with center $C'$ and
%radius $|C'Z|$,
passing through $Z$,
and let $O_{C'}$ intersect $O_{C}$ at $R$ and $R'$ (with $R$
%the point
on the same side as $Q$). This repetition used two more moves. Next we extend $BQ$ and $CR$ (now up to move 6) and let them meet at $S$. Finally, use the seventh move to draw the desired Apollonian circle $O_{S}$ centered at $S$
and passing through $Q$;
% with radius $%SP=
%|SQ|=|SR|$;
%, where $P$ is the point of tangency of $O_{S}$ with $O_{A}$;
see \figref{fig:4}.

%\begin{rmk}
%Note that it actually costs less to construct the fourth Apollonian circle than  the third one!
%\end{rmk}

\begin{rmk}
If a pair of lines, e.g.,
 $AC$ and $XZ$, are parallel (so $B'$ is at infinity), then  use the line $BY$ in lieu of $O_{B'}$ (the former is the limit of  the latter as $B'\to\infty$).
\end{rmk}

\subsection*{The Proof}  There is a unique circle such that inversion through it fixes $O_B$ and sends $O_A$ to $O_C$; we claim that $O_{B'}$ is this circle.  Indeed, such an inversion must send $X$ to $Z$, so its center must lie on $XZ$.  Its center also lies on the line perpendicular to $O_A$ and $O_C$, which is the line $AC$; thus its center is $B'=XZ\cap AC$.  Finally, the point $Y$ is fixed by this inversion, giving the claim.

Next it is easy to see  that the point of tangency of $O_{B}$ and the Apollonian circle $O_{S}$ must also lie on this inversion circle $O_{B'}$ (in which case  this point must be $Q$ as constructed). Indeed, since the inversion preserves the initial configuration of three circles, it must also fix  $O_{S}$, and hence also its point of tangency with $O_{B}$.

Finally, since $O_{B}$ and $O_{S}$ are tangent at $Q$, their centers are collinear with $Q$; that is, $S$ lies on the line $BQ$. The rest is elementary.\hfill $\sq$

\begin{rmk}\label{rmk:second}
%***Rewrite this***
The second solution $O_{S'}$ to the Apollonian problem can now be constructed in a further {\it three} moves. Indeed, %$O_{A'}$ intersects $O_{A}$ at two points, $P$ and $P'$,  say,  and similarly $O_{B'}$ intersects $O_{B}$ at $Q$ and $Q'$. T
the extra points of tangency $Q'$ and $R'$ are already on the page.
%cost us nothing.
Extend $BQ'$ and $CR'$ (two more moves); these intersect at $S'$, and drawing the circle $O_{S'}$ centered at $S'$ %having radius $|S'Q'|=|S'R'|$
and passing through $Q'$
costs a third move.
\end{rmk}

\begin{rmk}
Let $A'=BC\cap YZ$ be constructed similarly to $B'$ and $C'$.
Note that the triangles $\gD ABC$ and $\gD XYZ$ are perspective from the Gergonne point%
\footnote{See, e.g., Wikipedia for any (standard) terms not defined here and below.}%
.
By Desargue's theorem, they are therefore perspective from a line, which is $A'B'C'$, the so-called {\it Gergonne line}; see Oldknow \cite{Oldknow1996}, who seems to have been just shy of discovering the construction presented here.
\end{rmk}

\section{Other Constructions}

Apollonius's own solution did not survive antiquity \cite{Heath1981} and we only know of its existence through a ``mathscinet review'' by Pappus half a millennium later; perhaps we have simply rediscovered his work.
Vi\`ete's original
 solution  through inversion (see, e.g., \cite{Sarnak2011})  is logically extremely elegant but takes countless elementary moves.
There are many others but we highlight two in particular.

 \subsection*{Gergonne}
 Gergonne's own solution to the {\it general} Apollonian problem (that is, when the given circles are not necessarily tangent) is perhaps closest to ours (but of course the problem he is solving is more complicated).  He begins by constructing the radical circle $O_I$ for the initial circles $O_A$, $O_B$, and $O_C$, and identifies the six points $X$, $X'$, $Y$, $Y'$, $Z$, and $Z'$, where it intersects the three original circles.  Those points are taken in order around $O_I$, with $Y'$ and $Z$ on $O_A$, $Z'$ and $X$ on $O_B$, and $X'$ and $Y$ on $O_C$.  In our configuration, the radical circle is the incircle of triange  $ABC$ and $X=X'$, $Y=Y'$, and $Z=Z'$.

Every pair of circles can be thought of as being similar to each other via a dilation through a point.  In general, there are two such dilations.  This gives us six points of similarity, which lie on four lines, the four {\it lines of similitude}.  Each line generates a pair of tangent circles.  In our configuration, the point $B'$ is the center of the dilation that sends $O_A$ to $O_C$.  Since $O_A$ and $O_C$ are tangent, there is only one dilation, so we get only one line of similitude, the Gergonne line.

The radical circle of $O_B$, $O_I$, and a pair of tangent circles is centered on the line of similitude, so is where $XZ'$ intersects that line.  In our configuration, that gives us $B'$.  The radical circle is the one that intersects $O_I$ perpendicularly, so in our configuration it goes through $Y$.

\subsection*{Eppstein}

The previously simplest solution to our problem seems to have been that of Eppstein \cite{Eppstein2001, Eppstein2001a}, which used eleven elementary moves to draw $O_{S}$. His construction finds the tangency point $Q$ by first dropping the perpendicular to $AC$ through $B$, and then connecting a second line from $Y$ to one of the two points of intersection of this perpendicular with $O_{B}$. This second line intersects $O_{B}$ at $Q$ (or $Q'$, depending on the choice of intersection point). Note that constructing a perpendicular line is not an elementary operation, costing 3 moves. The second line {\it is} elementary, so Eppstein can construct $Q$ in 4 moves, then $R$ in 4 more, then two more lines $BQ$ and $CR$ to get the center $S$, and finally the circle $O_{S}$ in a total of 11 moves. To construct the other solution, $O_{S'}$, using his method, it would cost another five moves (as opposed to our three; see \rmkref{rmk:second}), since one needs to draw two more lines to produce $Q'$ and $R'$ (whereas our construction gives these as a byproduct).

\subsection*{Challenge:}
Construct (a generic configuration of) four mutually tangent circles in the plane using fewer than $15$ ($=1+2+5+7$) moves. Or prove (as we suspect) that this is impossible!

\bibliographystyle{alpha}

\bibliography{AKbibliog}

\end{document}